\pdfoutput=1
\RequirePackage{ifpdf}
\ifpdf 
\documentclass[pdftex]{sigma}
\else
\documentclass{sigma}
\fi

\numberwithin{equation}{section}

\newtheorem{thm}{Theorem}[section]

\newtheorem{prop}[thm]{Proposition}

\newtheorem{cor}[thm]{Corollary}

{\theoremstyle{definition}

\newtheorem{rema}[thm]{Remark}

}

\begin{document}


\renewcommand{\PaperNumber}{039}

\FirstPageHeading

\ShortArticleName{Some Remarks on Very-Well-Poised ${}_8\phi_7$ Series}

\ArticleName{Some Remarks on Very-Well-Poised $\boldsymbol{{}_8\phi_7}$ Series}

\Author{Jasper V.~STOKMAN}

\AuthorNameForHeading{J.V.~Stokman}
\Address{Korteweg--de Vries
Institute for Mathematics, University of Amsterdam,\\
Science Park 904, 1098 XH Amsterdam, The Netherlands}

\Email{\href{mailto:j.v.stokman@uva.nl}{j.v.stokman@uva.nl}}

\URLaddress{\url{http://staff.science.uva.nl/~jstokman/}}

\ArticleDates{Received April 05, 2012, in f\/inal form June 18, 2012; Published online June 27, 2012}

\Abstract{Nonpolynomial basic hypergeometric
eigenfunctions of the Askey--Wilson second order dif\/ference operator
are known to be expressible as very-well-poised ${}_8\phi_7$ series.
In this paper we use this fact to derive various basic hypergeometric
and theta function identities. We relate most of them
to identities from the existing literature on basic hypergeometric series.
This leads for example to a new derivation of a known
quadratic transformation formula for very-well-poised ${}_8\phi_7$ series.
We also provide a link to
Chalykh's theory on
(rank one, $\textup{BC}$ type) Baker--Akhiezer functions.}


\Classification{33D15; 33D45}

\section{Introduction}

\looseness=-1
In this paper we present derivations
of various basic hypergeometric and theta function identities
using the interpretation of very-well-poised ${}_8\phi_7$ series
as eigenfunctions of the Askey--Wilson se\-cond order dif\/ference operator~$\mathcal{D}$.
For instance, we reobtain a nonstandard type three term transformation formula
for very-well-poised ${}_8\phi_7$ series~\cite{GM}
as the connection formula for
the asympto\-ti\-cally free eigenfunctions of~$\mathcal{D}$, we investigate
the eigenfunctions of~$\mathcal{D}$
with trivial quantum mo\-nodromy and relate them to
(rank one, $\textup{BC}$ type) Baker--Akhiezer functions from~\mbox{\cite{Ch,ChE}},
we rederive the quadratic transformation formula \cite[(3.5.10)]{GR}
for very-well-poised ${}_8\phi_7$ series, and we obtain various theta
function identities
by translating symmetries of the Askey--Wilson function and of the
asymptotically free eigenfunction of~$\mathcal{D}$ in terms of the
associated normalized  $c$-functions.

The most general family of orthogonal polynomials satisfying
a second order $q$-dif\/ference equation is the family of Askey--Wilson
polynomials~\cite{AW}.
For our purposes it is convenient to view the associated
dif\/ference equations as
eigenvalue equations for the second order Askey--Wilson second order
dif\/ference operator~$\mathcal{D}$ already mentioned in the previous paragraph.
The operator~$\mathcal{D}$
depends, besides on the dif\/ference step-size and the deformation parameter~$q$, on four additional free parameters.

Nonpolynomial basic hypergeometric eigenfunctions of $\mathcal{D}$
have been subject of study in va\-rious papers
(see, e.g., \cite{HI,IR,KS,R,Rq,StDAHA,S}). Important examples are the
Askey--Wilson function $\mathcal{E}(\cdot,z)$
and the asymptotically free eigenfunction $\Phi(\cdot,z)$
(the corresponding
eigenvalue depends in an explicit way on $q^z+q^{-z}\in\mathbb{C}$). They
are def\/ined provided that $0<|q|<1$.
They are selfdual eigenfunctions (the role of the argument
and $z$ is interchangeable).
They naturally arise in harmonic analysis on the quantum $\textup{SU}(1,1)$
group
and in the study of the double af\/f\/ine Hecke algebra
of type $C^\vee C_1$ (the rank one Koornwinder case), see, e.g.,
\cite{KS2} and~\mbox{\cite{StDAHA,StSph}} respectively. {}From this representation
theoretic viewpoint~$\mathcal{E}(\cdot,z)$ plays the role of the spherical
function and~$\Phi(\cdot,z)$ the role of the Harish-Chandra
series.

Ruijsenaars' $R$-function \cite{R}
is another nonpolynomial selfdual eigenfunction of $\mathcal{D}$
which is required to satisfy yet another second order dif\/ference
equation of Askey--Wilson type.
The step-direction of the two Askey--Wilson second order
dif\/ference equations are allowed to be
co-linear (which corresponds to deformation parameter $q$ being
of modulus one).
The $R$-function arises as matrix coef\/f\/icient of representations of
the quantum double of the quantized universal enveloping
algebra of $\mathfrak{sl}_2$, see~\cite{vdB}.

The Askey--Wilson function $\mathcal{E}(\cdot,z)$ can be
explicitly expressed in terms of the asymptotically free eigenfunctions
$\Phi(\cdot,\pm z)$. The elliptic function $c(\cdot,z)$
governing the expansion coef\/f\/icients,
\[
\mathcal{E}(\cdot,z)=c(\cdot,z)\Phi(\cdot,z)+
c(\cdot,-z)\Phi(\cdot,-z),
\]
is called the (normalized) $c$-function \cite{KS}. It is explicitly given as
quotient of theta functions. The selfduality of $\mathcal{E}$
and $\Phi$ and the fact that $\mathcal{E}(-x,z)=
\mathcal{E}(x,z)$ then allow us to express $\Phi(-x,z)$
in terms of $\Phi(x,\pm z)$ (connection formula).
The cases when the connection coef\/f\/icient formula trivializes is particularly
interesting since it directly relates to the theory of Baker--Akhiezer
functions~\cite{Ch,ChE}. We discuss this in Section \ref{trivialsection}.

\looseness=-1
Suitable two parameter specializations of the Askey--Wilson polynomials
yield the continuous $q$-Jacobi polynomials. They have appeared
in two guises (see~\cite{Ra} and~\cite{AW}) which are interrelated by a
quadratic transformation formula going back to Singh~\cite{Si},
see also \cite[\S~4]{AW} and \cite[\S~3.10]{GR}. This quadratic
transformation formula was derived in \cite{AW}
using the orthogonality relations
of the Askey--Wilson polynomials.
Ruijsenaars~\cite{Rq} stressed that for
these parameter specializations the Askey--Wilson
second order dif\/ference operator $\mathcal{D}$
factorizes up to an additive constant as a~square of an Askey--Wilson
type second order dif\/ference operator with step-size half of
the step-size of $\mathcal{D}$. We use this observation
to prove a quadratic transformation
formula for a~two parameter family of
the asymptotically free eigenfunction $\Phi(\cdot,z)$ of $\mathcal{D}$.
This complements Ruijsenaars' results~\cite{Rq}, where he lifted
the quadratic transformation formula for continuous $q$-Jacobi polynomials
to a quadratic transformation formula for a two parameter
subfamily of the $R$-function.

Using the known explicit expression of $\Phi(\cdot,z)$
in terms of basic hypergeometric series we link the above mentioned
results
to various known identities for very-well-poised ${}_8\phi_7$ series.
For example, the quadratic transformation formula for $\Phi(\cdot,z)$
becomes the known quadratic transformation
formula \cite[(3.5.10)]{GR} for very-well-poised ${}_8\phi_7$ series
after applying suitable transformation formulas to both sides of the
identity, see Remark~\ref{Link}$(i)$ (this was observed by M.~Rahman).

Combining symmetries
of the Askey--Wilson function with its $c$-function expansion
yields nontrivial identities
for the $c$-function, hence nontrivial theta function identities.
For instance, the fact that $\mathcal{E}(\cdot,z)$
is invariant under negating the argument yields a theta function
identity which is a six parameter dependent
subcase
of a theta function identity \cite[Exercise~5.22]{GR}
due to Slater \cite{Sl}. The fact that $\Phi(\cdot,\pm z)$
satisf\/ies a quadratic transformation formula
but the Askey--Wilson function $\mathcal{E}(\cdot,z)$
does not, yields a nontrivial
identity for the $c$-function and consequently
a quadratic type theta function identity~\eqref{thetaident}.

The results in the present paper play
an important role in the study of the spectral problem of the
trigonometric
Macdonald--Ruijsenaars--Cherednik commuting family of dif\/ference
ope\-ra\-tors associated to root systems. 
They are for instance needed in the asymptotic analysis of
$q$-analogues of Harish-Chandra series
associated to root systems (cf.~\cite{LS,vM,vMS,StSph})
along codimension one facets of the Weyl chamber.
I will return to this topic in a future work.

\section[Eigenfunctions of the Askey-Wilson second order difference operator]{Eigenfunctions of the Askey--Wilson\\ second order dif\/ference operator}

We use, besides a deformation parameter
$0<q=e^{2\pi\sqrt{-1}\tau}<1$ ($\tau\in\sqrt{-1}\mathbb{R}_{>0}$)
and a choice of step-size $s\in\mathbb{Q}_{>0}$,
four free parameters $\{\kappa,\lambda,\upsilon,\varsigma\}$
which we call Hecke parameters (this name comes from their interpretation
as multiplicity parameters in the Cherednik--Macdonald theory of the
double af\/f\/ine Hecke algebra of type $C^\vee C_1$, see~\cite{NS}).
We assume that the
Hecke parameters $\kappa$, $\lambda$, $\upsilon$ and $\varsigma$ are real.
We set $q^x:=e^{2\pi\sqrt{-1}\tau x}$
for $x\in\mathbb{C}$.
In the paper $q$ will be f\/ixed throughout. Step-sizes
$s$ and $\frac{s}{2}$ will simultaneously appear in formulas. One can restrict
without loss of generality to considering step-sizes $2$ and $1$,
but formulas are more transparent when arbitrary values $s$ of the
step size
are taken into account because it makes the $s$-dependence
of the Askey--Wilson parameters
\[
\{a_s,b_s,c_s,d_s\}:=
\big\{q^{\kappa+\lambda},-q^{\kappa-\lambda},
q^{\frac{s}{2}+\upsilon+\varsigma},
-q^{\frac{s}{2}+\upsilon-\varsigma}\big\}
\]
explicit. We def\/ine dual Askey--Wilson parameters by
\[
\big\{\widetilde{a}_s,\widetilde{b}_s,\widetilde{c}_s,\widetilde{d}_s\big\}:=
\big\{q^{\kappa+\upsilon},-q^{\kappa-\upsilon},q^{\frac{s}{2}+\lambda+\varsigma},
-q^{\frac{s}{2}+\lambda-\varsigma}\big\},
\]
i.e.\ the roles of the Hecke parameters
$\lambda$ and $\upsilon$ are interchanged. The parameters
$a_s$, $b_s$, $\widetilde{a}_s$ and~$\widetilde{b}_s$ do not depend
on $s$, we therefore occasionally omit the subindex $s$
for these Askey--Wilson parameters.
Interchanging $\lambda$ and $\upsilon$ def\/ines an involution
on the set of Hecke parameters, hence also on the associated set of
Askey--Wilson parameters. In addition we have
\begin{gather*}
\widetilde{a}_s^2 =q^{-s}abcd,\qquad
\widetilde{a}_s\widetilde{b}_s =
a_sb_s,\qquad \frac{q^s\widetilde{a}_s}{\widetilde{b}_s}=c_sd_s,\\
\widetilde{a}_s\widetilde{c}_s =a_sc_s,
\qquad \frac{q^s\widetilde{a}_s}{\widetilde{c}_s}=b_sd_s,\qquad
\widetilde{a}_s\widetilde{d}_s =a_sd_s,\qquad
\frac{q^s\widetilde{a}_s}{\widetilde{d}_s}=b_sc_s.
\end{gather*}
We now f\/irst recall the
asymptotically free eigenfunction of the Askey--Wilson \cite{AW}
second order $q^s$-dif\/ference operator, which we will regard
here as a second order dif\/ference operator with step size $s$.
Explicitly, the Askey--Wilson second order dif\/ference operator $\mathcal{D}$,
acting on meromorphic functions on $\mathbb{C}$, is def\/ined by
\begin{gather*}
\bigl(\mathcal{D}f\bigr)(x) :=A(x)\bigl(f(x+s)-f(x)\bigr)+
A(-x)\bigl(f(x-s)-f(x)\bigr),\\
A(x) :=\frac{(1-a_sq^x)(1-b_sq^x)(1-c_sq^x)(1-d_sq^x)}{\widetilde{a}_s(1-q^{2x})
(1-q^{s+2x})}.
\end{gather*}
Sometimes it is important to write explicitly the dependence on the
parameters, in which case we write $\mathcal{D}$ as
$\mathcal{D}_{\kappa,\lambda,\upsilon,\varsigma;q}^{(s)}$.

\begin{rema}
The Askey--Wilson second order dif\/ference operator
$\mathcal{D}$ can be interpreted as
a second order $q^s$-dif\/ference operator
when acting on
$\tau^{-1}$-translation invariant
meromorphic functions on $\mathbb{C}$ (which are the meromorphic functions
of the form $g(q^x)$ with $g$ meromorphic on~$\mathbb{C}^*$).
\end{rema}

We now def\/ine the elementary function
$W(x,z)=W(x,z;\kappa,\lambda,\upsilon,\varsigma;q,s)$ by
\begin{gather*}
W(x,z):=q^{(\kappa+\lambda+x)(\kappa+\upsilon+z)/s}.
\end{gather*}
Note that $W(x+s,z)=q^{\kappa+\upsilon+z}W(x,z)=\widetilde{a}q^zW(x,z)$.

For generic $b_j$
the ${}_{r+1}\phi_r$ basic hypergeometric series is def\/ined by
the convergent power series
\[
{}_{r+1}\phi_r\left(
\begin{matrix} a_1,a_2,\ldots,a_{r+1}\\ b_1,b_2,\ldots,b_r\end{matrix};
q,z\right):=
\sum_{j=0}^{\infty}\frac{\bigl(a_1,a_2,\ldots,a_{r+1};q\bigr)_j}
{\bigl(q,b_1,\ldots,b_r;q\bigr)_j}z^j,\qquad |z|<1,
\]
where $\bigl(a_1,\ldots,a_s;q\bigr)_j=
\prod\limits_{r=1}^s\prod\limits_{i=0}^{j-1}(1-a_rq^i)$ for $j\in\mathbb{Z}_{\geq 0}\cup
\{\infty\}$ (empty products are equal to one by convention).
The very-well-poised ${}_8\phi_7$ series is def\/ined by
\begin{gather*}
{}_8W_7\bigl(\alpha_0;\alpha_1,\alpha_2,\alpha_3,\alpha_4,\alpha_5;q,z\bigr) =
{}_8\phi_7\left(\begin{matrix}
\alpha_0,q\alpha_0^{\frac{1}{2}},-q\alpha_0^{\frac{1}{2}},
\alpha_1,\ldots,\alpha_5\\
\alpha_0^{\frac{1}{2}},-\alpha_0^{-\frac{1}{2}},q\alpha_0/\alpha_1,\ldots,
q\alpha_0/\alpha_5\end{matrix};q,z\right)\\
\hphantom{{}_8W_7\bigl(\alpha_0;\alpha_1,\alpha_2,\alpha_3,\alpha_4,\alpha_5;q,z\bigr)}{}
=\sum_{r=0}^{\infty}\frac{1-\alpha_0q^{2r}}{1-\alpha_0}z^r
\prod_{j=0}^5\frac{\bigl(\alpha_j;q\bigr)_r}
{\bigl(q\alpha_0/\alpha_j;q\bigr)_r}.
\end{gather*}
In case $z=\alpha_0^2q^2/\alpha_1\alpha_2\alpha_3\alpha_4\alpha_5$
it has a meromorphic continuation to $\bigl(\mathbb{C}^*\bigr)^6$
as function of $(\alpha_0,\alpha_1,\ldots,\alpha_5)$. This follows from
the identity \cite[(III.36)]{GR} expressing such very-well-poised
${}_8\phi_7$ series as a sum of two ${}_4\phi_3$ series.
Def\/ine the holomorphic function $\textup{St}(x)=
\textup{St}(x;\kappa,\lambda,\upsilon,\varsigma;q,s)$ (``$\textup{St}$''
is standing for singular term) in $x\in\mathbb{C}$ by
\begin{gather*}
\textup{St}(x):=\bigl(q^{s+x}/a_s,q^{s+x}/b_s,q^{s+x}/c_s,q^{s+x}/d_s;
q^s\bigr)_{\infty}.
\end{gather*}
We write $\textup{St}^d(z):=\textup{St}(z;\kappa,\upsilon,\lambda,\varsigma;
q,s)$ for the singular term with respect to dual parameters,
\[
\textup{St}^d(z)=\bigl(q^{s+z}/\widetilde{a}_s,
q^{s+z}/\widetilde{b}_s,q^{s+z}/\widetilde{c}_s,q^{s+z}/\widetilde{d}_s;
q^s\bigr)_{\infty}.
\]
The following proposition combines and ref\/ines
observations from \cite{IR,KS,StSph}.

\begin{prop}\label{AWcase}
There exist unique
holomorphic functions $\Gamma_r$ on $\mathbb{C}$
$(r\geq 0)$ satisfying the following three
conditions,
\begin{enumerate}\itemsep=0pt
\item[$(1)$] $\Gamma_0(z)=\bigl(q^{s+2z};q^s\bigr)_{\infty}$.
\item[$(2)$] The power series $\Psi(x,z):=\sum\limits_{r=0}^{\infty}
\Gamma_r(z)q^{rx}$ is normally convergent on compacta of
$(x,z)\in\mathbb{C}\times\mathbb{C}$ $($consequently
 $\Psi(x,z)$ is a holomorphic
function in $(x,z)\in\mathbb{C}\times\mathbb{C})$.
\item[$(3)$] The meromorphic function
\begin{gather}\label{explicitPhi}
\Phi(x,z):=\frac{W(x,z)}{\textup{St}(x)\textup{St}^d(z)}\Psi(x,z)
\end{gather}
satisfies
\begin{gather}\label{equationAW}
\bigl(\bigl(\mathcal{D}_x-q^z-q^{-z}+\widetilde{a}+\widetilde{a}^{-1}\bigr)
\Phi\bigr)(x,z)=0,
\end{gather}
where $\mathcal{D}_x$ stands for the Askey--Wilson second order
difference operator
$\mathcal{D}=\mathcal{D}_{\kappa,\lambda,\upsilon,\varsigma;q}^{(s)}$
acting on the x-variable.
\end{enumerate}
Furthermore, $\Gamma_r$ is $\tau^{-1}$-translation invariant
and
\begin{gather}
\Psi(x,z) =\frac{\bigl(\frac{q^{s+x+z}a_s}{\widetilde{a}_s},
\frac{q^{s+x+z}b_s}{\widetilde{a}_s},\frac{q^{s+x+z}c_s}{\widetilde{a}_s},
\frac{q^{s+x+z}\widetilde{a}_s}{d_s},q^{s+2z},d_sq^x;q^s\bigr)_{\infty}}
{\bigl(\frac{q^{2s+x+2z}}{d_s};q^s\bigr)_{\infty}}\nonumber\\
\hphantom{\Psi(x,z) =}{}
\times{}_8W_7\left(\frac{q^{s+x+2z}}{d_s};\frac{q^{s+z}}{\widetilde{a}_s},
\frac{q^{s+z}}{\widetilde{d}_s},\widetilde{b}_sq^z,\widetilde{c}_sq^z,
\frac{q^{s+x}}{d_s};q^s,d_sq^x\right)\label{Psi8phi7}
\end{gather}
if $|d_sq^x|<1$.
\end{prop}

\begin{proof}
The explicit expression \eqref{Psi8phi7} is
\begin{gather}
 \Psi(x,z)=\frac{\bigl(\frac{q^{s+x+z}a_s}{\widetilde{a}_s},
\frac{q^{s+x+z}b_s}{\widetilde{a}_s},\frac{q^{s+x+z}c_s}{\widetilde{a}_s},
\frac{q^{s+x+z}\widetilde{a}_s}{d_s},q^{s+2z},d_sq^x;q^s\bigr)_{\infty}}
{\bigl(\frac{q^{s+x+2z}}{d_s};q^s\bigr)_{\infty}}\nonumber\\
\phantom{\Psi(x,z)=}{}
\times\sum_{r=0}^{\infty}
\left(1-\frac{q^{s+2sr+x+2z}}{d_s}\right)
\frac{\bigl(\frac{q^{s+z}}{\widetilde{a}_s},
\frac{q^{s+z}}{\widetilde{d}_s},\widetilde{b}_sq^z,
\widetilde{c}_sq^z,\frac{q^{s+x}}{d_s};q^s\bigr)_rd_s^rq^{rx}}
{\bigl(\frac{q^{s+x+z}\widetilde{a}_s}{d_s},
\frac{q^{s+x+z}a_s}{\widetilde{a}_s},
\frac{q^{s+x+z}b_s}{\widetilde{a}_s},
\frac{q^{s+x+z}c_s}{\widetilde{a}_s},q^{s+2z};q^s\bigr)_r}.\!\!\label{asymptoticform}
\end{gather}
It is a well def\/ined meromorphic function in
$(x,z)\in\mathbb{C}\times\mathbb{C}$ provided that $|d_sq^x|<1$,
with possible poles at $q^{s+sr+x+2z}=d_s$ ($r\in\mathbb{Z}_{\geq 0}$).
It can be expressed as sum of two
${}_4\phi_3$ series using \cite[(III.36)]{GR} with parameters
$(a,b,c,d,e,f,q)$ in \cite[(III.36)]{GR} specialized to
\[
(q^{s+x+2z}/d_s, q^{s+z}/\widetilde{a}_s,q^{s+z}/\widetilde{d}_s,
\widetilde{b}_sq^z,\widetilde{c}_sq^z,q^{s+x}/d_s,q^s),
\]
leading to the expression
\begin{gather}
\Psi(x,z) =\frac{\bigl(\frac{q^{s+x}}{a_s},d_sq^x,
\frac{q^{s+z}}{\widetilde{b}_s},
\frac{q^{s+z}}{\widetilde{c}_s},\frac{q^{s+x+z}a_s}{\widetilde{a}_s},
\frac{q^{s+x+z}\widetilde{a}_s}{d_s};q^s\bigr)_{\infty}}
{\bigl(\frac{d_s}{a_s};q^s\bigr)_{\infty}}
\nonumber\\
\hphantom{\Psi(x,z) =}{}
\times {}_4\phi_3\left(\begin{matrix}
a_sq^x,\frac{q^{s+x}}{d_s},\widetilde{b}_sq^z,\widetilde{c}_sq^z\\
\frac{q^{s+x+z}a_s}{\widetilde{a}_s},\frac{q^{s+x+z}\widetilde{a}_s}{d_s},
\frac{q^sa_s}{d_s}\end{matrix};
q^s,q^s\right)\nonumber\\
\hphantom{\Psi(x,z) =}{} +
\frac{\bigl(a_sq^x,\frac{q^{s+x}}{d_s},\widetilde{b}_sq^z,\widetilde{c}_sq^z,
\frac{q^{s+x+z}\widetilde{a}_s}{a_s},
\frac{q^{s+x+z}\widetilde{d}_s}{a_s};q^s\bigr)_{\infty}}
{\bigl(\frac{a_s}{d_s};q^s\bigr)_{\infty}}\nonumber\\
\hphantom{\Psi(x,z) =}{}
\times{}_4\phi_3\left(\begin{matrix}
\frac{q^{s+x}}{a_s},d_sq^x,\frac{q^{s+z}}{\widetilde{b}_s},
\frac{q^{s+z}}{\widetilde{c}_s}\\
\frac{q^{s+x+z}\widetilde{a}_s}{a_s},\frac{q^{s+x+z}\widetilde{d}_s}{a_s},
\frac{q^sd_s}{a_s}\end{matrix};q^s,q^s\right).\label{4phi3expression}
\end{gather}
This alternative expression provides the meromorphic continuation
and shows that $\Psi(x,z)$ is holomorphic in
$(x,z)\in\mathbb{C}\times\mathbb{C}$. The expression
\eqref{asymptoticform} of $\Psi(x,z)$
shows that $\Psi(x,z)$ satisf\/ies~(1) and~(2).
By~\cite{IR} (see also~\cite{KS} for notations close to the present one),
the resulting meromorphic function $\Phi(x,z)$ (see~\eqref{explicitPhi})
indeed satisf\/ies the dif\/ference equation~\eqref{equationAW}.

It remains to prove uniqueness. If a series of the form
\[
\Phi(x,z)=\frac{W(x,z)}{\textup{St}(x)\textup{St}^d(z)}\sum_{r\geq 0}^{\infty}
\Gamma_r(z)q^{rx}
\]
is a formal solution of \eqref{equationAW} then
the $\Gamma_r(z)$ ($r\geq 0$) satisfy recursion relations of the form
\[
\widetilde{a}(q^{sr}-1)(q^z-q^{-sr-z})\Gamma_r(z)=
\sum_{t=0}^{r-1}v_t^r(q^z)\Gamma_t(z),\qquad r\geq 1
\]
for some Laurent polynomials $v_t^r$. This shows that the $\Gamma_r(z)$
($r\geq 1$) are uniquely determined by $\Gamma_0(z)$.
\end{proof}

If confusion may arise about the parameter
dependencies then we write the asymptotically free solution
$\Phi(x,z)$ as
$\Phi(x,z;\kappa,\lambda,\upsilon,\varsigma;q,s)$ and
$\Psi(x,z)$ as $\Psi(x,z;\kappa,\lambda,\upsilon,\varsigma;q,s)$.

\begin{rema}\label{symmetric}
The characterization of $\Phi(x,z)$ as eigenfunction
of $\mathcal{D}$ is equivalent to
a characterization of $\Psi(x,z)=\sum\limits_{r=0}^{\infty}\Gamma_r(z)q^{rx}$
as eigenfunction of the second order dif\/ference operator obtained from
$\mathcal{D}$ by gauging it with gauge factor $W(\cdot,z)/\textup{St}(\cdot)$.
The gauged dif\/ference operator and the relevant
eigenvalue are symmetric under arbitrary permutations
of the Askey--Wilson parameters $a_s$, $b_s$, $c_s$, $d_s$. Since the
normalization $\Gamma_0(z)=\bigl(q^{s+2z};q^s\bigr)_{\infty}$
of $\Psi(x,z)$ is independent of the Askey--Wilson parameters,
it follows that $\Psi(x,z)$ is symmetric under arbitrary permutations
of the Askey--Wilson parameters $a_s$, $b_s$, $c_s$, $d_s$.
\end{rema}

Let $\theta(u;q):=\bigl(u,q/u;q\bigr)_{\infty}$ be the modif\/ied
Jacobi theta function and write
\[
\theta(u_1,\ldots,u_r;q)=
\prod\limits_{i=1}^r\theta(u_i;q)
\] for products of theta functions.
Def\/ine the (normalized) $c$-function
$c(x,z)=c(x,z;\kappa, \lambda,\upsilon,\varsigma;q,s)$ by
\begin{gather}\label{cfunction}
c(x,z):=\frac{\theta\bigl(\widetilde{a}_sq^{-z},\widetilde{b}_sq^{-z},
\widetilde{c}_sq^{-z},\frac{d_sq^{x-z}}{\widetilde{a}_s};q^s\bigr)}
{W(x,z)\theta\bigl(q^{-2z},d_sq^x;q^s\bigr)}.
\end{gather}
Using
$\theta(qu;q)=-u^{-1}\theta(u;q)$ it follows that
$c(x+s,z)=c(x,z)$ and $c(x,z+s)=c(x,z)$. We write $c^d(x,z)=
c(x,z;\kappa,\upsilon,\lambda,\varsigma;q,s)$ for
the $c$-function with respect to dual parameters.

The Askey--Wilson function $\mathcal{E}(x,z)=
\mathcal{E}(x,z;\kappa,\lambda,\upsilon,\varsigma;q,s)$
is def\/ined by
\begin{gather}
\mathcal{E}(x,z):=
 \frac{\bigl(\frac{\widetilde{a}_sq^{s+z-x}}{d_s},
\frac{\widetilde{a}_sq^{s+z+x}}{d_s},
a_sb_s,a_sc_s,\frac{q^sa_s}{d_s};q^s\bigr)_{\infty}}
{\bigl(\frac{q^{s+x}}{d_s},\frac{q^{s-x}}{d_s},\frac{q^{s+z}}{\widetilde{d}_s},
\widetilde{a}_s\widetilde{b}_s\widetilde{c}_sq^z;q^s\bigr)_{\infty}}\nonumber\\
 \hphantom{\mathcal{E}(x,z):=}{}
 \times{}_8W_7\left(\widetilde{a}_s\widetilde{b}_s\widetilde{c}_sq^{-s+z};
a_sq^x,a_sq^{-x},\widetilde{a}_sq^z,\widetilde{b}_sq^z,\widetilde{c}_sq^z;q^s,
\frac{q^{s-z}}{\widetilde{d}_s}\right)\label{AWfunction}
\end{gather}
for $|q^{s-z}/\widetilde{d}_s|<1$ (see \cite{KS}).
Using \cite[(III.36)]{GR} to express
$\mathcal{E}(x,z)$ as sum of two ${}_4\phi_3$ series it follows that
$\mathcal{E}(\cdot,\cdot)$ has a meromorphic extension to
$\mathbb{C}\times\mathbb{C}$.
We write $\mathcal{E}^d$ and $\Phi^d$ for the meromorphic functions
$\mathcal{E}$ and $\Phi$ with respect to dual parameters.
The following properties of $\mathcal{E}$ and $\Phi$ are known
from \cite{KS,StSph} (cf.\ also~\cite{HI,IR,S}).

\begin{prop}\label{cfunctionprop}\qquad\null
\begin{enumerate}\itemsep=0pt
\item[$(i)$] $\mathcal{E}(x,z)=\mathcal{E}^d(z,x)$ $($selfduality$)$.
\item[$(ii)$] $\mathcal{E}(-x,z)=\mathcal{E}(x,z)$
and $\mathcal{E}(x,-z)=\mathcal{E}(x,z)$.
\item[$(iii)$] $\Phi(x,z)=\Phi^d(z,x)$ $($selfduality$)$.
\item[$(iv)$] $\mathcal{E}(x,z)=c(x,z)\Phi(x,z)+c(x,-z)\Phi(x,-z)$ $(c$-function
expansion$)$.
\end{enumerate}
\end{prop}

\begin{proof}
$(i)$ This follows from the transformation formula
\cite[(III.23)]{GR} for very-well-poised ${}_8\phi_7$ series.

$(ii)$ By the explicit expression \eqref{AWfunction} it is clear
that $\mathcal{E}(-x,z)=\mathcal{E}(x,z)$. By
$(i)$ it then also follows that $\mathcal{E}(x,-z)=\mathcal{E}(x,z)$.

$(iii)$ This follows again by application of the transformation
formula \cite[(III.23)]{GR}
for very-well-poised ${}_8\phi_7$ series.

$(iv)$ Use Bailey's three term transformation formula
\cite[(III.37)]{GR} for very-well-poised ${}_8\phi_7$ series.
\end{proof}

\begin{rema}
The selfduality of $\Phi(x,z)$ and $\mathcal{E}(x,z)$
ensures that $\Phi(x,\cdot)$ and $\mathcal{E}(x,\cdot)$ are
eigenfunctions of the Askey--Wilson second order
dif\/ference operator with respect
to dual parameters.
\end{rema}

Note that the Askey--Wilson second order
dif\/ference operator $\mathcal{D}$ is invariant
under $x\mapsto -x$, hence $\Phi(-x,z)$ is again an eigenfunction
of $\mathcal{D}_x$ with eigenvalue $q^z+q^{-z}+\widetilde{a}+\widetilde{a}^{-1}$.
It can be expressed in terms of the eigenfunctions
$\Phi(x,z)$ and $\Phi(x,-z)$ as follows.

\begin{cor}
\begin{gather}\label{connectionformula}
\Phi(-x,z)=\left(\frac{c(x,z)-c^d(z,x)}{c^d(z,-x)}\right)\Phi(x,z)+
\frac{c(x,-z)}{c^d(z,-x)}\Phi(x,-z)
\end{gather}
$($connection formula$)$.
\end{cor}

\begin{proof}
By the previous proposition we have
\begin{gather*}
\mathcal{E}(x,z) =\mathcal{E}^d(z,-x)
 =c^d(z,-x)\Phi^d(z,-x)+c^d(z,x)\Phi^d(z,x)\\
\phantom{\mathcal{E}(x,z)}{} =c^d(z,x)\Phi(x,z)+c^d(z,-x)\Phi(-x,z).
\end{gather*}
Compared with the $c$-function expansion for $\mathcal{E}(x,z)$
(see Proposition \ref{cfunctionprop}$(iv)$)
we get the desired result.
\end{proof}

\begin{rema}
The connection formula \eqref{connectionformula} is not a
direct consequence of Bailey's
three term transformation formula \cite[(III.37)]{GR} for
very-well-poised ${}_8\phi_7$ series. This is
ref\/lected by the fact that the coef\/f\/icient of $\Phi(x,z)$ in~\eqref{connectionformula} does not admit an explicit
expression as a single product of theta functions. The connection
formula \eqref{connectionformula} is though directly related to
the three term transformation formula~\cite[(5.8)]{GM}.
\end{rema}

\begin{cor}\label{cquadratic}
The $c$-function satisfies
\begin{gather}
c(x,z)c^d(z,-x) c^d(-z,-x)+c(-x,z)c^d(z,x)c^d(-z,-x) \nonumber\\
\qquad{} =c(x,z)c(-x,z)c^d(-z,-x)+c(x,z)c(-x,-z)c^d(z,-x).\label{cqformula}
\end{gather}
\end{cor}

\begin{proof}
Since $\mathcal{E}(-x,z)=\mathcal{E}(x,z)$ we have
\[
c(x,z)\Phi(x,z)+c(x,-z)\Phi(x,-z)=
c(-x,z)\Phi(-x,z)+c(-x,-z)\Phi(-x,-z).
\]
Applying \eqref{connectionformula} twice to the right hand
side of this formula implies
\begin{gather}\label{twosides}
\alpha(x,z)\Phi(x,z)=-\alpha(x,-z)\Phi(x,-z)
\end{gather}
with the function $\alpha(x,z)$
given by
\[
\alpha(x,z)=\frac{c(-x,z)c(x,z)-c(-x,z)c^d(z,x)}{c^d(z,-x)}+
\frac{c(-x,-z)c(x,z)}{c^d(-z,-x)}-c(x,z).
\]
Replace in~\eqref{twosides} the variable $x$ by $x+ms$ and
consider the asymptotic behaviour as $m\rightarrow\infty$
of both sides, using the fact that~$c(x,z)$
is $s$-translation invariant in both $x$ and $z$ and using that
\[
\widetilde{a}^{-m}q^{-mz}\Phi(x+ms,z)=\Gamma_0(z)\bigl(1+\mathcal{O}(q^{sm})\bigr)
\]
as $m\rightarrow\infty$. It gives
$\alpha(x,z)=0$. This is equivalent to \eqref{cqformula}.
\end{proof}

\begin{rema}
Substituting in \eqref{cqformula} the explicit expression \eqref{cfunction}
of the $c$-function $c(x,z)$ gives the theta function identity
\begin{gather*}
 \theta\left(\frac{d}{\widetilde{a}}q^{x+z}\!,
\frac{d}{\widetilde{a}}q^{x-z},aq^x,bq^x,cq^x,dq^{-x},q^{2z};q\right)\!
 -
\theta\left(\frac{d}{\widetilde{a}}q^{-x+z},\frac{d}{\widetilde{a}}q^{x+z},
\widetilde{a}q^z,\widetilde{b}q^z,\widetilde{c}q^z,\widetilde{d}q^{-z},
q^{2x};q\right)\!\\
 \qquad{}=q^{2x}\theta\left(\frac{d}{\widetilde{a}}q^{-x-z},
\frac{d}{\widetilde{a}}q^{-x+z},aq^{-x},bq^{-x},cq^{-x},dq^x,q^{2z};q\right)\\
 \qquad\quad{} -
q^{2z}\theta\left(\frac{d}{\widetilde{a}}q^{-x-z},\frac{d}{\widetilde{a}}q^{x-z},
\widetilde{a}q^{-z},\widetilde{b}q^{-z},\widetilde{c}q^{-z},\widetilde{d}q^z,
q^{2x};q\right),
\end{gather*}
where we have taken $s=1$ and have written $a=a_1,\ldots,d=d_1$
(and similarly for the dual parameters). This is
a special case of Slater's theta function identity \cite[Exercise~5.22]{GR},
with the parameters $(a,b,c,d,e,f,g,h)$ in
\cite[Exercise~5.22]{GR} specialized to
\[
\left(1,\frac{d}{\widetilde{a}}q^{-x+z},\widetilde{a}q^{-z},
\widetilde{b}q^{-z},\widetilde{c}q^{-z},\frac{q^{1-z}}{\widetilde{d}},
\frac{d}{\widetilde{a}}q^{x+z},q^{2z}\right).
\]
Note that Slater's formula \cite[Exercise~5.22]{GR} is more general since it
has, besides $q$, seven free parameters,
while the formula~\eqref{cqformula} only has six.
\end{rema}

\section{The case of trivial quantum monodromy}\label{trivialsection}

Determining the monodromy representation for
Gauss' second order hypergeometric dif\/ferential equation
is equivalent to deriving its
connection coef\/f\/icient formulas. The connection coef\/f\/icient formulas
explicitly relate the fundamental series expansion solutions around the three
regular singularities of the hypergeometric dif\/ferential equation.
These formulas turn out to be directly related to well known three term
transformation formulas for the Gauss' hypergeometric function~${}_2F_1$.
From the above notion of monodromy only its incarnation in terms of
connection coef\/f\/icient formulas generalizes to the
dif\/ference setup (see Sauloy~\cite{Sau} and references therein for
a detailed discussion of this issue). We therefore say
that the explicit connection coef\/f\/icient formula~\eqref{connectionformula}
solves the quantum monodromy problem of~$\mathcal{D}$. In addition we
say that the quantum monodromy
is trivial if the coef\/f\/icient of~$\Phi(x,z)$ in~\eqref{connectionformula}
vanishes. The latter terminology is motivated as follows.
Firstly, for trivial quantum monodromy the connection coef\/f\/icient
formula \eqref{connectionformula} reduces to a simple equivariance
property of the following renormalization $\widetilde{\Phi}(x,z)
=\widetilde{\Phi}(x,z;\kappa,\lambda,\upsilon,\varsigma;
q,s)$ of the asymptotically free eigenfunction $\Phi(x,z)$,
\[
\widetilde{\Phi}(x,z):=c(x,z)\Phi(x,z),
\]
see Proposition \ref{trivial}$(iii)$ below. Secondly, for an important
subclass with trivial quantum monodromy, a suitable renormalization
of the asymptotically free eigenfunction $\Phi(x,z)$ has a terminating series
expansion, see Remark~\ref{BArem}$(i)$ for further discussions on
this issue.

Note that the renormalization $\widetilde{\Phi}(x,z)$ of
$\Phi(x,z)$ still satisf\/ies the eigenvalue equation
\[
\mathcal{D}\bigl(\widetilde{\Phi}(\cdot,z)\bigr)=
\big(q^z+q^{-z}-\widetilde{a}-\widetilde{a}^{-1}\big)\widetilde{\Phi}(\cdot,z)
\]
since $c(x,z)$ is $s$-translation invariant in $x$.
Write $\widetilde{\Phi}^d(z,x):=\widetilde{\Phi}(z,x;\kappa,\upsilon,\lambda,
\varsigma;q,s)$.

\begin{prop}\label{trivial}\quad\null
\begin{enumerate}\itemsep=0pt
\item[$(i)$] If $\frac{\kappa}{s},\frac{\lambda}{s},\frac{\upsilon}{s},
\frac{\varsigma}{s}\in
\frac{1}{2}\mathbb{Z}$ with an even number of them
being integers, then $c(x,z)=c^d(z,x)$.
\end{enumerate}
In the remaining items of the proposition
we assume that $(\kappa,\lambda,\upsilon,\varsigma)$
is a four-tuple of real parameters such that $c(x,z)=c^d(z,x)$.
\begin{enumerate}\itemsep=0pt
\item[$(ii)$] $\widetilde{\Phi}(x,z)=\widetilde{\Phi}^d(z,x)$
$($selfduality$)$.
\item[$(iii)$] $\widetilde{\Phi}(-x,z)=\widetilde{\Phi}(x,-z)$.
\item[$(iv)$] $\mathcal{E}(x,z)=\widetilde{\Phi}(x,z)+
\widetilde{\Phi}(-x,z)$.
\end{enumerate}
\end{prop}

\begin{proof}
$(i)$ Consider the quotient $\alpha:=c(x,z)/c^d(z,x)$
(we suppress the dependence on the variables and parameters).
Using the explicit
expression of the normalized $c$-function~\eqref{cfunction}
in terms of theta functions it follows that $\alpha=1$ if
$\bigl(\frac{\kappa}{s},\frac{\lambda}{s},
\frac{\upsilon}{s},\frac{\varsigma}{s})$ is taken from the set
\begin{gather*}
 \left\{\bigl(0,0,0,0\bigr), \left(\frac{1}{2},\frac{1}{2},0,0\right),
\left(\frac{1}{2},0,\frac{1}{2},0\right),
\left(\frac{1}{2},0,0,\frac{1}{2}\right),\right.\\
\left. \qquad\qquad\qquad \left(0,\frac{1}{2},\frac{1}{2},0\right),
\left(0,\frac{1}{2},0,\frac{1}{2}\right),
\left(0,0,\frac{1}{2},\frac{1}{2}\right),
\left(\frac{1}{2},\frac{1}{2},\frac{1}{2},
\frac{1}{2}\right)\right\}.
\end{gather*}
In addition, it is easy to check that
$\alpha$ is invariant under integral shifts of the rescaled
Hecke parameters
$\frac{\kappa}{s}$, $\frac{\lambda}{s}$, $\frac{\upsilon}{s}$,
$\frac{\varsigma}{s}$. This gives the result.

$(ii)$ This follows from Proposition \ref{cfunctionprop}$(iii)$.

$(iii)$  This follows from the connection formula
\eqref{connectionformula} for $\Phi(x,z)$.

$(iv)$ This is immediate from the $c$-function expansion of
the Askey--Wilson function, see Proposition \ref{cfunctionprop}$(iv)$.
\end{proof}

\begin{rema}\label{BArem}\qquad\null
\begin{enumerate}\itemsep=0pt
\item[$(i)$] Due to \cite[\S~4]{LS}
the function $\widetilde{\Phi}(x,z)$ for Hecke parameters
$(\kappa,\lambda,\upsilon,\varsigma)$ satisfying,
besides the integrality conditions from Proposition~\ref{trivial}$(i)$,
suitable additional (positivity) conditions, is up to normalization
Chalykh's $\textup{BC}_1$ type normalized
Baker--Akhiezer function~\cite[\S~6]{Ch}. An important property of the Baker--Akhiezer
function is the fact that it has a terminating series expansion. Surprisingly
this is not evident from the explicit expressions of $\widetilde{\Phi}(x,z)$
in terms of basic hypergeometric series. The properties
from Proposition~\ref{trivial} can
be easily matched with the properties (see \cite{Ch,ChE})
of the Baker--Akhiezer functions.
For instance, the Askey--Wilson function $\mathcal{E}(x,z)$
relates to the symmetrized normalized Baker--Akhiezer function
$\Phi_+$ from \cite[(3.20)]{ChE}. In particular,
formula \cite[Theorem~3.9]{ChE} relating the symmetrized Baker--Akhiezer function
to Askey--Wilson polynomials matches with the known
fact that the Askey--Wilson function $\mathcal{E}(x,z)$
reduces to the normalized Askey--Wilson polynomial for suitable discrete
values of $z$, see \eqref{AWfctpol} below.
\item[$(ii)$] The relation between the $\textup{A}_1$ type Baker--Akhiezer
function (see \cite[\S~4.1]{Ch})
and Heine's basic hypergeometric series ${}_2\phi_1$ was stressed
in an informal note of Koornwinder \cite{Knote}.
\end{enumerate}
\end{rema}

\section[The factorization of the Askey-Wilson second order
difference operator]{The factorization of the Askey--Wilson\\ second order
dif\/ference operator}

Ruijsenaars \cite{Rq}
analyzed when the square
$\bigl(\mathcal{D}+\widetilde{a}+\widetilde{a}^{-1}\bigl)^2$
is an Askey--Wilson second order dif\/ference operator again (with doubled
step-size). An important special case turns out to be
when the Hecke parameters are of the form
$(\kappa,\lambda,0,0)$.
In our notations the resulting
formula \cite[(3.11)]{Rq} reads as
\begin{gather}
\bigl(\mathcal{D}_{\kappa,\lambda,0,0;q}^{(s/2)}-q^{\frac{z}{2}}-q^{-\frac{z}{2}}+
q^\kappa+q^{-\kappa}\bigr)
 \bigl(\mathcal{D}_{\kappa,\lambda,0,0;q}^{(s/2)}+
q^{\frac{z}{2}}+q^{-\frac{z}{2}}+
q^\kappa+q^{-\kappa}\bigr)\nonumber\\
 \qquad{}
=\mathcal{D}_{\kappa,\lambda,\kappa,\lambda;q}^{(s)}-q^z-q^{-z}+q^{2\kappa}+
q^{-2\kappa}.\label{factorizationinitial}
\end{gather}

\begin{rema}\label{relatienotatie}
Our operator $\mathcal{D}+\widetilde{a}+\widetilde{a}^{-1}=
\mathcal{D}+q^{2\kappa}+q^{-2\kappa}$ for arbitrary Hecke parameters
$(\kappa,\lambda,\upsilon,\varsigma)$
is essentially the second order dif\/ference operator
$A_\delta(\mathbf{c};\cdot)$ in Ruijsenaars' paper~\cite[(2.1)]{Rq}
with the parameters $(c_0,c_1,c_2,c_3)$ in \cite{Rq}
related to our Hecke parameters
by $\kappa+\lambda=-\sqrt{-1}c_0$, $\kappa-\lambda=-\sqrt{-1}c_1$,
$\upsilon+\varsigma=-\sqrt{-1}c_2$ and $\upsilon-\varsigma=-\sqrt{-1}c_3$.
\end{rema}

A simple derivation of \eqref{factorizationinitial} is as follows.
Consider the second order dif\/ference operator
$\mathcal{L}=\mathcal{L}_{\kappa,\lambda;q}^{(s/2)}$ def\/ined by
\[
\bigl(\mathcal{L}f\bigr)(x):=\frac{(1-q^{\kappa+\lambda+x})
(1+q^{\kappa-\lambda+x})}{q^\kappa(1-q^{2x})}
f\left(x+\frac{s}{2}\right)+\frac{(1-q^{\kappa+\lambda-x})
(1+q^{\kappa-\lambda-x})}{q^\kappa(1-q^{-2x})}
f\left(x-\frac{s}{2}\right).
\]
Then
\begin{gather*}
\bigl(\bigl(\mathcal{L}_{\kappa,\lambda,q}^{(s/2)}-
q^\kappa-q^{-\kappa}\bigr)f\bigr)(x) =
\frac{(1-q^{\kappa+\lambda+x})(1+q^{\kappa-\lambda+x})}{q^\kappa(1-q^{2x})}
\left(f\left(x+\frac{s}{2}\right)-f(x)\right)\\
 \hphantom{\bigl(\bigl(\mathcal{L}_{\kappa,\lambda,q}^{(s/2)}-
q^\kappa-q^{-\kappa}\bigr)f\bigr)(x) =}{}
 +\frac{(1-q^{\kappa+\lambda-x})(1+q^{\kappa-\lambda-x})}{q^\kappa(1-q^{-2x})}
\left(f\left(x-\frac{s}{2}\right)-f(x)\right)\\
 \hphantom{\bigl(\bigl(\mathcal{L}_{\kappa,\lambda,q}^{(s/2)}-
q^\kappa-q^{-\kappa}\bigr)f\bigr)(x)}{}
 =\bigl(\mathcal{D}_{\kappa,\lambda,0,0;q}^{(s/2)}f\bigr)(x).
\end{gather*}
Formula \eqref{factorizationinitial}
thus is equivalent to
\begin{gather}\label{factorization}
\bigl(\mathcal{L}_{\kappa,\lambda;q}^{(s/2)}-q^{\frac{z}{2}}-q^{-\frac{z}{2}}\bigr)
\bigl(\mathcal{L}_{\kappa,\lambda;q}^{(s/2)}+q^{\frac{z}{2}}+q^{-\frac{z}{2}}\bigr)
=\mathcal{D}_{\kappa,\lambda,\kappa,\lambda;q}^{(s)}
-q^z-q^{-z}+q^{2\kappa}+q^{-2\kappa},
\end{gather}
which is an easy check (use that the Askey--Wilson parameters
associated to Hecke parameters $(\kappa,\lambda,\kappa,\lambda)$
satisfy $c_s=q^{\frac{s}{2}}a_s=
q^{\frac{s}{2}+\kappa+\lambda}$ and
$d_s=q^{\frac{s}{2}}b_s=-q^{\frac{s}{2}+\kappa-\lambda}$).

The eigenfunctions of $\mathcal{D}_x$ of the form $p_n(q^x+q^{-x})$
($n\in\mathbb{Z}_{\geq 0}$) with~$p_n$ a polynomial of degree~$n$,
are the well known Askey--Wilson polynomials~\cite{AW} (see also Section~\ref{polred}). In the special case that the associated four Hecke parameters
are taken to be $(\kappa,\lambda,\kappa,\lambda)$ (corresponding to the
condition that the Askey--Wilson parameters satisfy
$c_s=q^{\frac{s}{2}}a_s$ and $d_s=q^{\frac{s}{2}}b_s$) the Askey--Wilson
polynomials are called the continuous $q$-Jacobi polynomials
by Askey and Wilson \cite[\S~4]{AW}.
In \cite[(4.22)]{AW} a quadratic transformation formula
for balanced ${}_4\phi_3$ series, going back to Singh~\cite{Si},
is used to relate the above continuous $q$-Jacobi polynomials
to Rahman's def\/inition~\cite{Ra} of the continuous $q$-Jacobi polynomials
(see also Section~\ref{polred}).

Ruijsenaars \cite{Rq}
used the factorization~\eqref{factorizationinitial}
to motivate and analyze quadratic transformation
formulas for the hyperbolic nonpolynomial
generalization of the continuous $q$-Jacobi polynomial,
which is Ruijsenaars' $R$-function with the continuous $q$-Jacobi
specialization $(\kappa,\lambda,\kappa,\lambda)$
of the associated Hecke algebra parameters (see
Remark \ref{relatienotatie}
for the relation with Ruijsenaars' notations~\cite{Rq}).
In the following section we
use the factorization \eqref{factorizationinitial} to prove quadratic
transformation formulas for the asymptotically free
eigenfunction $\Phi(\cdot,z;\kappa,\lambda,\kappa,\lambda;q,s)$.

\section{Quadratic transformation formulas}

For a function $\Xi(x,z;\kappa,\lambda,\upsilon,\varsigma;q,s)$ we write
\[
\Xi_R(x,z):=\Xi\left(x,\frac{z}{2};\kappa,\lambda,0,0;q,\frac{s}{2}\right).
\]
This is the parameter specialization which reduces the Askey--Wilson
polynomials to Rahman's version of the
continuous $q$-Jacobi polynomials~\cite{Ra} (in base $q^{\frac{s}{2}}$).
Furthermore we set $\Xi_R^d(x,z):=
\Xi(x,\frac{z}{2};\kappa,0,\lambda,0;q,\frac{s}{2})$ for its dual version.
Since
\begin{gather*}
\mathcal{L}_{\kappa,\lambda,q}^{(s/2)}=\mathcal{D}_{\kappa,\lambda,0,0;q}^{(s/2)}
+q^\kappa+q^{-\kappa},
\end{gather*}
we know from the previous section that the meromorphic function
$\Phi_R(x,z)$ satisf\/ies
\begin{gather}\label{differenceBC}
\bigl(\mathcal{L}_x-q^{\frac{z}{2}}-q^{-\frac{z}{2}}\bigr)\Phi_R(\cdot,z)=0.
\end{gather}
In fact, $\Phi_R(x,z)$ is the unique solution to
\eqref{differenceBC}
which is of the form
\[
\Phi_R(x,z)=\frac{W_R(x,z)}{\textup{St}_R(x)\textup{St}_R^d(z)}
\sum_{r=0}^{\infty}\Gamma_{R,r}(z)q^{rx}
\]
with $\Gamma_{R,r}(z)$ holomorphic in $z\in\mathbb{C}$, with
$\Gamma_{R,0}(z)=\bigl(q^{\frac{s}{2}+z};q^{\frac{s}{2}}\bigr)_{\infty}$
and with the power series $\Psi_R(x,z):=\sum\limits_{r=0}^{\infty}\Gamma_{R,r}(z)q^{rx}$
converging normally on compacta of $(x,z)\in\mathbb{C}\times\mathbb{C}$.

For a function $\Xi(x,z;\kappa,\lambda,\upsilon,\varsigma;q,s)$
we write
\[
\Xi_J(x,z):=\Xi(x,z;\kappa,\lambda,\kappa,\lambda;q,s).
\]
This is the parameter specialization which reduces the Askey--Wilson
polynomials to Askey's and Wilson's version of the
continuous $q$-Jacobi polynomials~\cite[\S~4]{AW}.
In addition we write
$\Xi_J^d(x,z):=\Xi(x,z;\kappa,\kappa,\lambda,\lambda;q,s)$ for its
dual version.

\begin{thm}\label{QTasymptotic}
$\Phi_J(x,z)=\Phi_R(x,z)$.
\end{thm}

\begin{proof}
First of all note that the Askey--Wilson parameters associated to
the Hecke parameters $(\kappa,\lambda,0,0)$, deformation parameter $q$
and step-size $\frac{s}{2}$ are given by
\[\big(q^{\kappa+\lambda},-q^{\kappa-\lambda},
q^{\frac{s}{4}},-q^{\frac{s}{4}}\big).
\]
The expression of $\Psi(x,z)$ as a sum of two ${}_4\phi_3$'s,
see \eqref{4phi3expression}, shows that
\[
\bigl(-q^{\frac{s}{4}+x};q^{\frac{s}{2}}\bigr)_{\infty}^{-1}\Psi_R(x,z)
\]
is holomorphic in $(x,z)\in\mathbb{C}\times\mathbb{C}$. By Remark
\ref{symmetric} also
\[
\bigl(q^{\frac{s}{4}+x},-q^{\frac{s}{4}+x};q^{\frac{s}{2}}\bigr)_{\infty}^{-1}
\Psi_R(x,z)
=\bigl(q^{\frac{s}{2}+2x};q^s\bigr)_{\infty}^{-1}\Psi_R(x,z)
\]
is holomorphic in $(x,z)\in\mathbb{C}\times\mathbb{C}$. We write its power
series expansion in $q^x$ as
\[
\bigl(q^{\frac{s}{2}+2x};q^s\bigr)_{\infty}^{-1}\Psi_R(x,z)=
\sum_{r=0}^{\infty}\widetilde{H}_r(z)q^{rx}.
\]
The $\widetilde{H}_r(z)$ are holomorphic in $z\in\mathbb{C}$, the series
converges normally on compacta of $(x,z)\in\mathbb{C}\times\mathbb{C}$,
and $\widetilde{H}_0(z)=\bigl(q^{\frac{s}{2}+z};q^{\frac{s}{2}}\bigr)_{\infty}$.
Def\/ine new holomorphic functions by
\[
H_r(z):=\bigl(-q^{\frac{s}{2}+z};q^{\frac{s}{2}}\bigr)_{\infty}\widetilde{H}_r(z)
\]
for $r\in\mathbb{Z}_{\geq 0}$. Then $H_0(z)=\bigl(q^{s+2z};q^s\bigr)_{\infty}$
and
\[
\Phi_R(x,z)=\frac{W_R(x,z)\bigl(q^{\frac{s}{2}+2x};q^s\bigr)_{\infty}}
{\textup{St}_R(x)\textup{St}_R^d(z)
\bigl(-q^{\frac{s}{2}+z};q^{\frac{s}{2}}\bigr)_{\infty}}
\sum_{r=0}^{\infty}H_r(z)q^{rx}
\]
with the series converging normally on compacta of $(x,z)\in\mathbb{C}\times
\mathbb{C}$. A direct computation now shows that
\[
\Phi_R(x,z)=\frac{W_J(x,z)}{\textup{St}_J(x)
\textup{St}_J^d(z)}\sum_{r=0}^{\infty}
H_r(z)q^{rx}.
\]
By \eqref{differenceBC} and \eqref{factorization} we furthermore have
\[
\bigl(\mathcal{D}_{\kappa,\lambda,\kappa,\lambda;q}^{(s)}-q^z-q^{-z}
+q^{2\kappa}+q^{-2\kappa}\bigr)\Phi_R(\cdot,z)=0.
\]
Hence $\Phi_R(x,z)$ satisf\/ies the characterizing properties
of $\Phi_J(x,z)$.
\end{proof}

By the proof of the theorem we in particular have
\[
\Psi_R(x,z)=\frac{\bigl(q^{\frac{s}{2}+2x};q^s\bigr)_{\infty}}
{\bigl(-q^{\frac{s}{2}+z};q^{\frac{s}{2}}\bigr)_{\infty}}\Psi_J(x,z).
\]
Substituting the explicit expressions of $\Psi_R$ and $\Psi_J$
in this formula
gives the following quadratic transformation formula for very-well-poised
${}_8\phi_7$ series.
\begin{cor}
\begin{gather}
 \frac{\bigl(q\beta xz,-\frac{qxz}{\beta},
\frac{q^{\frac{3}{2}}xz}{\alpha},-q^{\frac{1}{2}}\alpha xz;q\bigr)_{\infty}}
{\bigl(q^{\frac{1}{2}}x,-q^{\frac{3}{2}}xz^2,
\frac{q^2\beta xz^2}{\alpha};q\bigr)_{\infty}}
\,
  {}_8W_7\left(-q^{\frac{1}{2}}xz^2;\frac{qz}{\alpha},
-\frac{q^{\frac{1}{2}}z}{\beta},-\alpha z, q^{\frac{1}{2}}\beta z,
-q^{\frac{1}{2}}x;q,-q^{\frac{1}{2}}x\right)\nonumber\\
\qquad{} =\frac{\bigl(-\frac{q^2xz^2}{\alpha\beta},-q\alpha\beta xz^2,
-\frac{q\alpha x}{\beta};q^2\bigr)_{\infty}}
{\bigl(-\frac{q^3\beta xz^4}{\alpha};q^2\bigr)_{\infty}}\nonumber\\
 \qquad\quad{} \times {}_8W_7\left(
-\frac{q\beta xz^4}{\alpha};\frac{q^2z^2}{\alpha^2},
-qz^2,-z^2,q\beta^2z^2,
-\frac{q\beta x}{\alpha};q^2,-\frac{q\alpha x}{\beta}
\right)\label{qtrans8W7}
\end{gather}
if both $|q^{\frac{1}{2}}x|<1$ and $|\frac{q\alpha x}{\beta}|<1$.
\end{cor}

\begin{rema}\label{Link}\qquad\null
\begin{enumerate}\itemsep=0pt
\item[$(i)$] As observed by Mizan Rahman (private communication),
the quadratic transformation formula \eqref{qtrans8W7}
for ${}_8W_7$ is equivalent to the quadratic transformation formula
\cite[(3.5.10)]{GR} by applying to the left hand side of~\eqref{qtrans8W7}
the transformation formula~\cite[(III.23)]{GR} with parameters
$(a,b,c,d,e,f)$ in \cite[(III.23)]{GR} specialized to
\[
\left(-q^{\frac{1}{2}}xz^2, -\frac{q^{\frac{1}{2}}z}{\beta},
-\alpha z, -q^{\frac{1}{2}}x,\frac{qz}{\alpha},q^{\frac{1}{2}}\beta z\right)
\]
and to the right hand side of \eqref{qtrans8W7} the transformation formula
\cite[(III.23)]{GR} with parameters $(a,b,c,d,e,f,q)$ in \cite[(III.23)]{GR}
specialized to
\[
\left(-\frac{q\beta xz^4}{\alpha},-qz^2,-z^2,q\beta^2z^2,
\frac{q^2z^2}{\alpha^2},-\frac{q\beta x}{\alpha},q^2\right).
\]
\item[$(ii)$] The formal classical limit $q\uparrow 1$ of \eqref{qtrans8W7}
can be computed after replacing the parameters
$(x,z,\alpha,\beta)$ in~\eqref{qtrans8W7}
by $(-x,q^z,q^\alpha,q^\beta)$ and moving all the inf\/inite $q$-shifted
factorials in \eqref{qtrans8W7} to one side.
The resulting classical limit turns out to be trivial, since
both sides reduce to
\[
{}_2F_1\left(\begin{matrix} 1+2z-2\alpha,\frac{1}{2}+2z+2\beta\\
1+4z\end{matrix}; \frac{4x}{(1+x)^2}\right),
\]
where ${}_2F_1$ is Gauss' hypergeometric series
(in contrast to the classical limits of the quadra\-tic transformations
of very-well-poised ${}_8\phi_7$'s from
\cite[\S~5]{RV}, which reduce to nontrivial quadratic
transformations for ${}_2F_1$).
On the polynomial level \eqref{qtrans8W7} reduces to the
quadratic transformation formula for the continuous $q$-Jacobi polynomials,
see Section~\ref{polred}, which is known to reduce to non\-tri\-vial quadratic
transformations (see~\cite[(4.24)]{AW}) on the classical level when taking
the limit $q\uparrow 1$ after replacing the parameters $(x,z,\alpha,\beta)$
by $(-q^x,q^z,q^\alpha,q^\beta)$.
\end{enumerate}
\end{rema}

There is a dual version of \eqref{qtrans8W7}, which can be obtained
by using the selfduality of $\Phi(x,z)$ to both sides of the quadratic
transformation formula $\Phi_J(x,z)=\Phi_R(x,z)$ before substituting the
explicit expression as a ${}_8W_7$ series. It leads to the following quadratic
transformation formula.
\begin{cor}
\begin{gather}
 \frac{\bigl(-q^{\frac{1}{2}}\alpha xz,-q^{\frac{1}{2}}\beta z,
\frac{q^{\frac{3}{2}}xz}{\alpha};q\bigr)_{\infty}}
{\bigl(-\frac{q^{\frac{3}{2}}x^2z}{\beta};q\bigr)_{\infty}}
\,
{}_8W_7\left(
-\frac{q^{\frac{1}{2}}x^2z}{\beta};\frac{qx}{\alpha\beta},
-q^{\frac{1}{2}}x,-\frac{\alpha x}{\beta},q^{\frac{1}{2}}x,
-\frac{q^{\frac{1}{2}}z}{\beta};q,-q^{\frac{1}{2}}\beta z\right)\nonumber\\
\qquad{} =\frac{\bigl(-q\alpha\beta xz^2, \frac{q^2\alpha xz^2}{\beta},
-\frac{q^2xz^2}{\alpha\beta},\frac{q^{3}\beta xz^2}{\alpha};q^2
\bigr)_{\infty}}
{\bigl(-q^3x^2z^2,-q^2z^2,\frac{q^2x^2z^2}{\beta^2};q^2\bigr)_{\infty}}\nonumber\\
 \qquad\quad{} \times
{}_8W_7\left(
-qx^2z^2; \frac{q^2x}{\alpha\beta},
-\frac{q\beta x}{\alpha},-\frac{\alpha x}{\beta},
q\alpha\beta x, -\frac{q\beta z^2}{\alpha};
q^2,-qx\right).\label{qtrans8W7twee}
\end{gather}
\end{cor}
Also \eqref{qtrans8W7twee} can be related to \cite[(3.5.10)]{GR}
by applying the transformation formula \cite[(III.23)]{GR}
to both sides, cf.\ Remark \ref{Link}$(i)$.

\section{Polynomial reduction}\label{polred}

Both the asymptotically free
eigenfunctions $\Phi(x,z)$ of the Askey--Wilson second order dif\/ference
operator, and the Askey--Wilson function $\mathcal{E}(x,z)$,
reduce to the Askey--Wilson polynomials when~$z$ is specialized
appropriately (see, e.g., \cite{KS,StSph}).
Concretely, for $\Phi(x,z)$ we have for $n\in\mathbb{Z}_{\geq 0}$,
\begin{gather}\label{AWpolynomialexpression}
\Phi(x,-\kappa-\upsilon-ns)=
\frac{a_s^{-2n}\bigl(a_sb_s,a_sc_s,a_sd_s;q^s\bigr)_n
\bigl(\frac{q^{2(1-n)s}}{a_sb_sc_sd_s};q^s\bigr)_{\infty}}
{\bigl(q^{(n-1)s}a_sb_sc_sd_s;q^s\bigr)_n\textup{St}^d(-\kappa-\upsilon-ns)}
P_n(x)
\end{gather}
with $P_n(x)=P_n(x;\kappa,\lambda,\upsilon,\varsigma;q,s)$ the
normalized Askey--Wilson polynomial in $q^x+q^{-x}$ of deg\-ree~$n$,
\begin{gather*}
P_n(x):={}_4\phi_3\left(
\begin{matrix} q^{-ns},q^{(n-1)s}a_sb_sc_sd_s,a_sq^x,a_sq^{-x}\\
a_sb_s,a_sc_s,a_sd_s\end{matrix};q^s,q^s\right).
\end{gather*}
This can be proved directly from the expression of
$\Phi(x,-\kappa-\upsilon-ns)$ as
a very-well-poised ${}_8W_7$ series (see Proposition~\ref{AWcase})
by f\/irst applying Watson's transformation \cite[(III.17)]{GR}
with parameters $(a,b,c,d,e,f,q)$ in~\cite[(III.17)]{GR} specialized to
\[
\left(\frac{q^{2(1-n)s+x}}{a_sb_sc_sd_s^2},
\frac{q^{s+x}}{d_s}, \frac{q^{(2-n)s}}{a_sb_sc_sd_s},
\frac{q^{(1-n)s}}{a_sd_s}, \frac{q^{(1-n)s}}{b_sd_s},
\frac{q^{(1-n)s}}{c_sd_s},q^s\right),
\]
followed by Sear's transformation \cite[(III.16)]{GR} with parameters
$(a,b,c,d,e,f,q)$ in \cite[(III.16)]{GR} specialized to
\[
\left(\frac{q^{(1-n)s}}{a_sd_s},\frac{q^{(1-n)s}}{b_sd_s},
\frac{q^{(1-n)s}}{c_sd_s},\frac{q^{2(1-n)s}}{a_sb_sc_sd_s},
\frac{q^{(1-n)s+x}}{d_s},\frac{q^{(1-n)s-x}}{d_s},q^s\right).
\]
Formula \eqref{AWpolynomialexpression}
can also be proved by observing that $P_n(x)$ is an
eigenfunction of $\mathcal{D}_x$ with eigenvalue
$\widetilde{a}(q^n+1)+\widetilde{a}^{-1}(q^{-n}+1)$ (see~\cite{AW})
of the form
\[
P_n(x)=\frac{\bigl(q^{-ns},q^{(n-1)s}a_sb_sc_sd_s;q^s\bigr)_n(-a_s)^n
q^{\frac{s}{2}n(n+1)}}
{\bigl(a_sb_s,a_sc_s,a_sd_s;q^s\bigr)_n}q^{nx}+\sum_{m<n}c_mq^{mx},
\]
hence up to a multiplicative constant it must be equal to
$\Phi(x,-\kappa-\upsilon-ns)$ (this argument generalizes to the setting
of Macdonald--Koornwinder polynomials, see~\cite{LS,StSph}).
Also the Askey--Wilson function $\mathcal{E}(x,-\kappa-\upsilon-ns)$
($n\in\mathbb{Z}_{\geq 0}$) is a multiple of the Askey--Wilson polynomial,
\begin{gather}\label{AWfctpol}
\mathcal{E}(x,-\kappa-\upsilon-ns)=
\frac{\bigl(a_sb_s,a_sc_s;q^s\bigr)_{\infty}}
{\bigl(\frac{q^s}{a_sd_s};q^s\bigr)_{\infty}}
P_n(x).
\end{gather}
This can again be proved directly using transformation formulas
and the expression of $\mathcal{E}(x,z)$
as a ${}_8W_7$ series, see~\cite{KS}. It can also be proved using
the $c$-function expansion of $\mathcal{E}(x,z)$
(see Proposition~\ref{cfunctionprop}$(iv)$)
and \eqref{AWpolynomialexpression}, since $c(x;\kappa+\upsilon+ns)=0$
for $n\in\mathbb{Z}_{\geq 0}$.

Specializing now $z=-2\kappa-ns$ ($n\in\mathbb{Z}_{\geq 0}$)
in the identity $\Phi_J(x,z)=\Phi_R(x,z)$ and using~\eqref{AWpolynomialexpression}
twice gives, after straightforward simplif\/ications,
the following result (we take $s=2$ without loss of generality),
\begin{gather*}
{}_4\phi_3\left(\begin{matrix}
q^{-2n},aq^x,aq^{-x},q^{2n}a^2b^2\\
ab,qa^2,qab
\end{matrix};q^2,q^2\right)=
{}_4\phi_3\left(\begin{matrix}
q^{-n},aq^x,aq^{-x},-q^{n}ab\\
ab,q^{\frac{1}{2}}a,-q^{\frac{1}{2}}a
\end{matrix};q,q\right)
\end{gather*}
for $n\in\mathbb{Z}_+$. This quadratic transformation formula
was f\/irst proved by Singh~\cite{Si}. It was reobtained by Askey and Wilson
\cite[(4.22)]{AW} and interpreted as quadratic transformation formula
\cite[(4.22)]{AW} for the continuous $q$-Jacobi polynomial (to get it in
the same form one has to apply Sear's transformation formula
\cite[(III.15)]{GR}), see also \cite[\S~3.10]{GR} for further discussions.
It was also obtained in \cite[(3.20)]{Rq}
($-a^2q$ in \cite[(3.20)]{Rq} should be $-abq$) by specialization of the
quadratic transformation formula \cite[(3.16)]{Rq} for the $R$-function.
A similar polynomial reduction can be done with
the dual version~\eqref{qtrans8W7twee} of the quadratic transformation
formula, in which case it reduces to the
quadratic transformation formula~\cite[(3.1)]{AW}.

\section{A theta function identity}

By the $c$-function
expansion of the Askey--Wilson function (see Proposition
\ref{cfunctionprop}$(iv)$) and by Theorem \ref{QTasymptotic} we have
\begin{gather*}
\mathcal{E}_J(x,z) =c_J(x,z)\Phi_J(x,z)+c_J(x,-z)\Phi_J(x,-z)\nonumber\\
\phantom{\mathcal{E}_J(x,z)}{}
=c_J(x,z)\Phi_R(x,z)+c_J(x,-z)\Phi_R(x,-z)
\end{gather*}
while on the other hand,
\[
\mathcal{E}_R(x,z)=c_R(x,z)\Phi_R(x,z)+c_R(x,-z)\Phi_R(x,-z).
\]
The $c$-functions $c_J(x,z)$ and $c_R(x,z)$ are explicitly given by
\begin{gather}
c_J(x,z) =\frac{\theta\bigl(q^{2\kappa-z},-q^{-z},q^{\frac{s}{2}+2\lambda-z},
-q^{\frac{s}{2}-\kappa-\lambda+x-z};q^s\bigr)}
{W_J(x,z)\theta\bigl(q^{-2z},-q^{\frac{s}{2}+\kappa-\lambda+x};q^s\bigr)},\nonumber\\
c_R(x,z) =\frac{\theta\bigl(q^{2\kappa-z};q^s\bigr)\theta\bigl(
q^{\frac{s}{4}+\lambda-\frac{z}{2}},-q^{\frac{s}{4}-\kappa+x-\frac{z}{2}};
q^{\frac{s}{2}}\bigr)}
{W_R(x,z)\theta\bigl(q^{-z},-q^{\frac{s}{4}+x};q^{\frac{s}{2}}\bigr)}.\label{cQR}
\end{gather}
Note that $W_J(x,z)=q^{(\kappa+\lambda+x)(2\kappa+z)/s}=W_R(x,z)$ and that
\[
\frac{c_J(x,z)}{c_R(x,z)}=\frac{
\theta\bigl(-q^{\frac{s}{2}-\kappa-\lambda+x-z};q^s\bigr)
\theta\bigl(-q^{\frac{s}{4}+\lambda-\frac{z}{2}},
-q^{\frac{3}{4}+x};q^{\frac{s}{2}}\bigr)}
{\theta\bigl(-q^{\frac{s}{2}-z},-q^{\frac{s}{2}+\kappa-\lambda+x};
q^s\bigr)\theta\bigl(-q^{\frac{s}{4}-\kappa+x-\frac{z}{2}};q^{\frac{s}{2}}
\bigr)}
\]
is not invariant under $z\rightarrow -z$ for generic Hecke parameters,
since it is zero
at $z=2\lambda+\frac{s}{2}+2\pi\sqrt{-1}\tau$ but nonzero at
$z=-2\lambda-\frac{s}{2}-2\pi\sqrt{-1}\tau$.
Hence $\mathcal{E}_J(x,z)$ cannot be
of the form $\alpha(x,z)\mathcal{E}_R(x,z)$
with $\alpha(x,z)$ meromorphic and $s$-invariant in the variable $x$.

\begin{rema}\label{not}
Both $\mathcal{E}_J(\cdot,z)$ and $\mathcal{E}_R(\cdot,z)$ are eigenfunctions
of $\mathcal{D}_{\kappa,\lambda,\kappa,\lambda;q}^{(s)}$ with eigenvalue
$q^z+q^{-z}-q^{2\kappa}-q^{-2\kappa}$ (for $\mathcal{E}_R(\cdot,z)$ this follows
from \eqref{factorizationinitial}
and the fact that $c_R(\cdot,z)$ is
$\frac{s}{2}$-translation invariant).
On the other hand, $\mathcal{E}_R(\cdot,z)$
is an eigenfunction of $\mathcal{L}_{\kappa,\lambda,q}^{(s/2)}$
with eigenvalue $q^{\frac{z}{2}}+q^{-\frac{z}{2}}$, but this is not true
for $\mathcal{E}_J(\cdot,z)$ since $c_J(\cdot,z)$ is not
$\frac{s}{2}$-translation invariant.
\end{rema}

\begin{rema}
Ruijsenaars' $R$-function $\mathcal{R}(\cdot,z)$ \cite{Rq}
is a hyperbolic eigenfunction of the
Askey--Wilson second order dif\/ference operator
(implying in particular that it
admits an analytic continuation
to the regime $|q|=1$) satisfying $\mathcal{R}(-x,z)=\mathcal{R}(x,z)$.
For the $R$-function $\mathcal{R}$ a quadratic
transformation formula of the form $\mathcal{R}_J(x,z)=\mathcal{R}_R(x,z)$
holds true, see \cite[(3.16)]{Rq}. The discrepancy with the fact
that $\mathcal{E}_J(x,z)\not=\mathcal{E}_R(x,z)$ is not unexpected
since the $R$-function $\mathcal{R}$ in the trigonometric regime $|q|<1$
has a nontrivial factorization in Askey--Wilson functions:
it expands as a sum of two terms, each term being essentially the product
of two Askey--Wilson functions, see \cite[Theorem~6.5]{vdBRS}.
\end{rema}

The quadratic transformation formula $\Phi_J(x,z)=\Phi_R(x,z)$
implies the following result for the connection coef\/f\/icients
in~\eqref{connectionformula}.

\begin{prop}
We have
\begin{gather}
\frac{c_J(x,-z)}{c_J^d(z,-x)} =\frac{c_R(x,-z)}{c_R^d(\frac{z}{2},-2x)},\qquad
\frac{c_J(x,z)-c_J^d(z,x)}{c_J^d(z,-x)} =
\frac{c_R(x,z)-c_R^d(\frac{z}{2},2x)}{c_R^d(\frac{z}{2},-2x)}.\label{quadratictheta}
\end{gather}
\end{prop}

\begin{proof}
Using \eqref{connectionformula} and Theorem \ref{QTasymptotic}
we have
\begin{gather}
 \left(\frac{c_J(x,z)-c_J^d(z,x)}{c_J^d(z,-x)}\right)\Phi_J(x,z)+
\frac{c_J(x,-z)}{c_J^d(z,-x)}\Phi_J(x,-z) \nonumber\\
 \qquad
=\left(\frac{c_R(x,z)-c_R^d(\frac{z}{2},2x)}{c_R^d(\frac{z}{2},-2x)}\right)
\Phi_J(x,z)+\frac{c_R(x,-z)}{c_R^d(\frac{z}{2},-2x)}\Phi_J(x,-z),\label{twoways}
\end{gather}
since both sides are equal to $\Phi_J(-x,z)=\Phi_R(-x,z)$. By a straightforward
a\-symptotic argument (compare with the proof of Corollary~\ref{cquadratic})
it follows that the coef\/f\/icients of $\Phi_J(x,z)$
(resp.\ of~$\Phi_J(x,-z)$) on both sides of this equation should be the
same.

Alternatively, one verif\/ies by a direct computation that
\[
\frac{c_J(x,-z)}{c_J^d(z,-x)}=
q^{-\frac{4\kappa x}{s}}q^{\frac{2(\kappa+\lambda)z}{s}}
\frac{\theta\bigl(q^{2x},q^{\frac{s}{2}+2\lambda+z},
q^{2\kappa+z};q^s\bigr)}
{\theta\bigl(q^z,q^{\kappa+\lambda+x},-q^{\kappa-\lambda+x};q^{\frac{s}{2}}\bigr)}
=\frac{c_R(x,-z)}{c_R^d(\frac{z}{2},-2x)},
\]
yielding the f\/irst equality of \eqref{quadratictheta}.
Combined with \eqref{twoways} this implies the second equality \linebreak of~\eqref{quadratictheta}.
\end{proof}

Since $c(x,z)$ is $s$-translation invariant in both $x$ and $z$,
it follows from the right hand sides of~\eqref{quadratictheta}
that the quotients
\[
\frac{c_J(x,-z)}{c_J^d(z,-x)},\qquad
\frac{c_J(x,z)-c_J^d(z,x)}{c_J^d(z,-x)}
\]
are $\frac{s}{2}$-translation invariant in $x$ (although $c_J(x,z)$ is
not, cf.\ Remark~\ref{not}).

Using \eqref{cQR} and using the explicit expressions
\begin{gather*}
c_J^d(z,x) =\frac{\theta\bigl(q^{\kappa+\lambda-x},
-q^{\kappa-\lambda-x},q^{\frac{s}{2}+\kappa+\lambda-x},
-q^{\frac{s}{2}-\kappa-\lambda+z-x};q^s\bigr)}
{W_J^d(z,x)\theta\bigl(q^{-2x},-q^{\frac{s}{2}+z};q^s\bigr)},\\
c_R^d(z,x) =\frac{\theta\bigl(q^{\kappa+\lambda-\frac{x}{2}},
-q^{\kappa-\lambda-\frac{x}{2}},q^{\frac{s}{4}-\frac{x}{2}},
-q^{\frac{s}{4}-\kappa+z-\frac{x}{2}};q^{\frac{s}{2}}\bigr)}
{W_R^d(z,x)\theta\bigl(q^{-x},-q^{\frac{s}{4}+\lambda+z};q^{\frac{s}{2}}\bigr)}
\end{gather*}
with $W_J^d(z,x)=q^{(\kappa+\lambda+x)(2\kappa+z)/s}$
and $W_R^d(z,x)=q^{2(\kappa+z)(\kappa+\lambda+\frac{x}{2})/s}$,
we have explicitly,
\begin{gather*}
\frac{c_J(x,z)-c_J^d(z,x)}{c_J^d(z,-x)} =
q^{-2x(2\kappa+z)/s}\frac{\theta\bigl(q^{2x};q^s\bigr)}
{\theta\bigl(q^{\kappa+\lambda+x},-q^{\kappa-\lambda+x},
q^{\frac{s}{2}+\kappa+\lambda+x},-q^{\frac{s}{2}-\kappa-\lambda+x+z};q^s\bigr)}\\
\hphantom{\frac{c_J(x,z)-c_J^d(z,x)}{c_J^d(z,-x)} =}{} \times\left\{
\frac{\theta\bigl(q^{2\kappa-z},-q^{-z},q^{\frac{s}{2}+2\lambda-z},
-q^{\frac{s}{2}-\kappa-\lambda+x-z},-q^{\frac{s}{2}+z};q^s\bigr)}
{\theta\bigl(q^{-2z},-q^{\frac{s}{2}+\kappa-\lambda+x};q^s\bigr)}\right.\\
 \left.\qquad
\hphantom{\frac{c_J(x,z)-c_J^d(z,x)}{c_J^d(z,-x)} =}{}
-\frac{\theta\bigl(q^{\kappa+\lambda-x},-q^{\kappa-\lambda-x},
q^{\frac{s}{2}+\kappa+\lambda-x},-q^{\frac{s}{2}-\kappa-\lambda-x+z};q^s\bigr)}
{\theta\bigl(q^{-2x};q^s\bigr)}\right\}
\end{gather*}
and
\begin{gather*}
\frac{c_R(x,z)-c_R^d(\frac{z}{2},2x)}{c_R^d(\frac{z}{2},-2x)} =
q^{-2x(2\kappa+z)/s}
\frac{\theta\bigl(q^{2x};q^{\frac{s}{2}}\bigr)}
{\theta\bigl(q^{\kappa+\lambda+x},-q^{\kappa-\lambda+x},
q^{\frac{s}{4}+x},-q^{\frac{s}{4}-\kappa+x+\frac{z}{2}};q^{\frac{s}{2}}\bigr)}\\
\hphantom{\frac{c_R(x,z)-c_R^d(\frac{z}{2},2x)}{c_R^d(\frac{z}{2},-2x)} =}{}
 \times\left\{
\frac{\theta\bigl(q^{2\kappa-z};q^s\bigr)
\theta\bigl(q^{\frac{s}{4}+\lambda-\frac{z}{2}},-q^{\frac{s}{4}-\kappa+x-\frac{z}{2}},
-q^{\frac{s}{4}+\lambda+\frac{z}{2}};q^{\frac{s}{2}}\bigr)}
{\theta\bigl(q^{-z},-q^{\frac{s}{4}+x};q^{\frac{s}{2}}\bigr)}\right.\\
 \left.\qquad
\hphantom{\frac{c_R(x,z)-c_R^d(\frac{z}{2},2x)}{c_R^d(\frac{z}{2},-2x)} =}{}
-\frac{\theta\bigl(q^{\kappa+\lambda-x},-q^{\kappa-\lambda-x},
q^{\frac{s}{4}-x},-q^{\frac{s}{4}-\kappa-x+\frac{z}{2}};q^{\frac{s}{2}}\bigr)}
{\theta\bigl(q^{-2x};q^{\frac{s}{2}}\bigr)}\right\}.
\end{gather*}
Hence \eqref{quadratictheta} is equivalent to
the following theta function identity.

\begin{cor}
\begin{gather}
\theta\left(a^2,b^2c^2,-b^2,qb^2d^2,-\frac{qb^2}{acd},
-qb^2;q^2\right)
 -\theta\left(b^4,-\frac{qc}{ad},acd,-\frac{ac}{d},
qacd,-\frac{qa}{b^2cd};q^2\right)\nonumber \\
 \qquad{} =\frac{\theta\bigl(-qab^2cd;q^2\bigr)\theta\bigl( -b^2;q\bigr)}
 {\theta\bigl(qa^2;q^2\bigr)\theta\bigl(-q^{\frac{1}{2}}abc; q\bigr)}
 \left\{
\theta\left(bc,-bc,a^2,q^{\frac{1}{2}}bd,-\frac{q^{\frac{1}{2}}b}{ac},
-\frac{q^{\frac{1}{2}}d}{b};q\right)\right.\nonumber\\
 \left.\qquad
\quad{}
 -
\theta\left(b^2,-q^{\frac{1}{2}}a,acd,-\frac{ac}{d},
q^{\frac{1}{2}}a,-\frac{q^{\frac{1}{2}}a}{bc};q\right)\right\}.\label{thetaident}
\end{gather}
\end{cor}

\subsection*{Acknowledgements}
I thank Tom Koornwinder for drawing
my attention to the quadratic transformation formula for
continuous $q$-Jacobi polynomials. I thank Mizan Rahman for
pointing out to me how the quadratic transformations \eqref{qtrans8W7}
and \eqref{qtrans8W7twee} for very-well-poised ${}_8\phi_7$
series are related to the known quadratic transformation formula
\cite[(3.5.10)]{GR} (see Remark \ref{Link}$(i)$).

\pdfbookmark[1]{References}{ref}
\LastPageEnding


\begin{thebibliography}{99}
\footnotesize\itemsep=0pt

\bibitem{AW}
Askey R., Wilson J., Some basic hypergeometric orthogonal polynomials that
  generalize {J}acobi polynomials, \textit{Mem. Amer. Math. Soc.} \textbf{54}
  (1985), no.~319.

\bibitem{vdB}
van~de Bult F.J., Ruijsenaars' hypergeometric function and the modular double
  of {${\mathcal U}_q(\mathfrak{sl}_2({\mathbb C}))$}, \href{http://dx.doi.org/10.1016/j.aim.2005.05.023}{\textit{Adv. Math.}}
  \textbf{204} (2006), 539--571, \href{http://arxiv.org/abs/math.QA/0501405}{math.QA/0501405}.

\bibitem{vdBRS}
van~de Bult F.J., Rains E.M., Stokman J.V., Properties of generalized
  univariate hypergeometric functions, \href{http://dx.doi.org/10.1007/s00220-007-0289-0}{\textit{Comm. Math. Phys.}} \textbf{275}
  (2007), 37--95, \href{http://arxiv.org/abs/math.CA/0607250}{math.CA/0607250}.

\bibitem{Ch}
Chalykh O.A., Macdonald polynomials and algebraic integrability, \href{http://dx.doi.org/10.1006/aima.2001.2033}{\textit{Adv.
  Math.}} \textbf{166} (2002), 193--259, \href{http://arxiv.org/abs/math.QA/0212313}{math.QA/0212313}.

\bibitem{ChE}
Chalykh O.A., Etingof P., Orthogonality relations and Cherednik identities for
  multivariable Baker--Akhiezer functions, \href{http://arxiv.org/abs/1111.0515}{arXiv:1111.0515}.


\bibitem{GR}
Gasper G., Rahman M., Basic hypergeometric series, \href{http://dx.doi.org/10.1017/CBO9780511526251}{\textit{Encyclopedia of
  Mathematics and its Applications}}, Vol.~96, 2nd ed., Cambridge University
  Press, Cambridge, 2004.

\bibitem{GM}
Gupta D.P., Masson D.R., Contiguous relations, continued fractions and
  orthogonality, \href{http://dx.doi.org/10.1090/S0002-9947-98-01879-0}{\textit{Trans. Amer. Math. Soc.}} \textbf{350} (1998),
  769--808, \href{http://arxiv.org/abs/math.CA/9511218}{math.CA/9511218}.

\bibitem{HI}
Haine L., Iliev P., Askey--{W}ilson type functions with bound states,
  \href{http://dx.doi.org/10.1007/s11139-006-8478-6}{\textit{Ramanujan~J.}} \textbf{11} (2006), 285--329, \href{http://arxiv.org/abs/math.QA/0203136}{math.QA/0203136}.

\bibitem{IR}
Ismail M.E.H., Rahman M., The associated {A}skey--{W}ilson polynomials,
  \href{http://dx.doi.org/10.2307/2001881}{\textit{Trans. Amer. Math. Soc.}} \textbf{328} (1991), 201--237.

\bibitem{KS}
Koelink E., Stokman J.V., The {A}skey--{W}ilson function transform,
  \href{http://dx.doi.org/10.1155/S1073792801000575}{\textit{Int. Math. Res. Not.}} (2001), no.~22, 1203--1227, \href{http://arxiv.org/abs/math.CA/0004053}{math.CA/0004053}.

\bibitem{KS2}
Koelink E., Stokman J.V., Fourier transforms on the quantum {${\rm SU}(1,1)$}
  group, \href{http://dx.doi.org/10.2977/prims/1145477332}{\textit{Publ. Res. Inst. Math. Sci.}} \textbf{37} (2001), 621--715,
  \href{http://arxiv.org/abs/math.QA/9911163}{math.QA/9911163}.

\bibitem{Knote}
Koornwinder T., Comment on the paper ``Macdonald polynomials and algebraic
  integrability'' by O.A.~Chalykh, available at
  \url{http://staff.science.uva.nl/~thk/art/comment/ChalykhComment.pdf}.

\bibitem{LS}
Letzter G., Stokman J.V., Macdonald dif\/ference operators and {H}arish-{C}handra
  series, \href{http://dx.doi.org/10.1112/plms/pdm055}{\textit{Proc. Lond. Math. Soc.~(3)}} \textbf{97} (2008), 60--96,
  \href{http://arxiv.org/abs/math.QA/0701218}{math.QA/0701218}.

  \bibitem{vM}
van Meer M., Bispectral quantum {K}nizhnik--{Z}amolodchikov equations for
  arbitrary root systems, \href{http://dx.doi.org/10.1007/s00029-010-0039-6}{\textit{Selecta Math.~(N.S.)}} \textbf{17} (2011),
  183--221, \href{http://arxiv.org/abs/0912.3784}{arXiv:0912.3784}.

\bibitem{vMS}
van Meer M., Stokman J., Double af\/f\/ine {H}ecke algebras and bispectral quantum
  {K}nizhnik--{Z}amolodchikov equations, \href{http://dx.doi.org/10.1093/imrn/rnp165}{\textit{Int. Math. Res. Not.}}  (2010),
  no.~6,  969--1040, \href{http://arxiv.org/abs/0812.1005}{arXiv:0812.1005}.

\bibitem{NS}
Noumi M., Stokman J.V., Askey--{W}ilson polynomials: an af\/f\/ine {H}ecke algebra
  approach, in Laredo {L}ectures on {O}rthogonal {P}olynomials and {S}pecial
  {F}unctions, \emph{Adv. Theory Spec. Funct. Orthogonal Polynomials}, Nova Sci.
  Publ., Hauppauge, NY, 2004, 111--144, \href{http://arxiv.org/abs/math.QA/0001033}{math.QA/0001033}.

\bibitem{Ra}
Rahman M., The linearization of the product of continuous {$q$}-{J}acobi
  polynomials, \href{http://dx.doi.org/10.4153/CJM-1981-076-8}{\textit{Canad.~J. Math.}} \textbf{33} (1981), 961--987.

\bibitem{RV}
Rahman M., Verma A., Quadratic transformation formulas for basic hypergeometric
  series, \href{http://dx.doi.org/10.2307/2154269}{\textit{Trans. Amer. Math. Soc.}} \textbf{335} (1993), 277--302.

\bibitem{R}
Ruijsenaars S.N.M., A~generalized hypergeometric function satisfying four
  analytic dif\/ference equations of {A}skey--{W}ilson type, \href{http://dx.doi.org/10.1007/s002200050840}{\textit{Comm. Math.
  Phys.}} \textbf{206} (1999), 639--690.

\bibitem{Rq}
Ruijsenaars S.N.M., Quadratic transformations for a function that generalizes
  {${}_2F_1$} and the {A}skey--{W}ilson polynomials, \href{http://dx.doi.org/10.1007/s11139-006-0257-x}{\textit{Ramanujan~J.}}
  \textbf{13} (2007), 339--364.

\bibitem{Sau}
Sauloy J., Syst\`emes aux {$q$}-dif\/f\'erences singuliers r\'eguliers:
  classif\/ication, matrice de connexion et mo\-no\-dromie, \textit{Ann. Inst.
  Fourier (Grenoble)} \textbf{50} (2000), 1021--1071.

\bibitem{Si}
Singh V.N., The basic analogues of identities of the {C}ayley--{O}rr type,
  \href{http://dx.doi.org/10.1112/jlms/s1-34.1.15}{\textit{J.~London Math. Soc.}} \textbf{34} (1959), 15--22.

\bibitem{Sl}
Slater L.J., A note on equivalent product theorems, \textit{Math. Gaz.}
  \textbf{38} (1954), 127--128.

\bibitem{StDAHA}
Stokman J.V., An expansion formula for the {A}skey--{W}ilson function,
  \href{http://dx.doi.org/10.1006/jath.2001.3647}{\textit{J.~Approx. Theory}} \textbf{114} (2002), 308--342,
  \href{http://arxiv.org/abs/math.CA/0105093}{math.CA/0105093}.

\bibitem{StSph}
Stokman J.V., The $c$-function expansion of a basic hypergeometric function
  associated to root systems, \href{http://arxiv.org/abs/1109.0613}{arXiv:1109.0613}.

\bibitem{S}
Suslov S.K., Some orthogonal very-well-poised {${}_8\phi_7$}-functions that
  generalize {A}skey--{W}ilson polynomials, \href{http://dx.doi.org/10.1023/A:1011439924912}{\textit{Ramanujan~J.}} \textbf{5}
  (2001), 183--218, \href{http://arxiv.org/abs/math.CA/9707213}{math.CA/9707213}.

\end{thebibliography}
\end{document}